# Three beliefs that lend illusory legitimacy to Cantor's diagonal argument

Bhupinder Singh Anand

## 1. Cantor's diagonal argument, Gödel's proof, and Turing's Halting problem

Whatever other beliefs there may remain for considering Cantor's diagonal argument[1] as mathematically legitimate, there are three that, prima facie, lend it an illusory legitimacy; they need to be explicitly discounted appropriately.

The first - Cantor's diagonal argument defines a non-countable Dedekind real number[2]; the second - Gödel uses the argument to define a formally undecidable, but interpretively true, proposition; and the third - Turing uses the argument to define an uncomputable Dedekind real number.

## 2. Cantor's diagonal argument

In the first case, we may define any natural number, expressed in binary notation, and followed by a period and a non-terminating sequence of the integers 0 and 1, as a Cantorian real number.

Cantor's diagonal argument, then, considers any, given, 1-1 correspondence:

(*)     $n <=> C_n$

where $n$ ranges over the natural numbers, and $C_n$ is a Cantorian real number of the form $0.a_{n1}a_{n2}a_{n3}...a_{ni}...$, where $a_{ni}$ is either 0 or 1 for any given $i$. Cantor then defines the

---

Cantorian number $0.b_1b_2b_3...b_i...$, such that $b_i$ is not equal to $a_{ii}$ for any $i$. Clearly, $0.b_1b_2b_3...b_i...$ does not correspond to any natural number $n$ in the given 1-1 correspondence (*). Since any given natural number $n$ can obviously be corresponded to the Cantorian number where $n$ is expressed in binary notation, and followed by a period and a non-terminating series of 0's, Cantor concludes that the Cantorian real numbers are uncountable.

It is not at all obvious, however, that the following conclusion can be drawn from the above argument: Every Cantorian real number necessarily defines a Dedekind cut[3] in the rational numbers (or its equivalent Cauchy sequence[4] of rationals).

For instance, consider the non-terminating series of the odd natural numbers 1, 3, 5, 7, ... We can express this as the non-terminating binary series:

    1, 11, 101, 111, ... ,

and correspond to it the non-terminating sequence of rational numbers:

    0.1, 0.11, 0.101, 0.111, ... .

Clearly, this last series oscillates; it can be computed by a classical Turing machine; and it defines a Cantorian real number if we accept Turing's implicit postulation[5] that every "circle-free a-machine" uniquely defines a non-terminating sequence of the integers 0 and 1. However, equally clearly, it does not yield a Cauchy sequence of rational numbers that can be taken to define[6] a Dedekind real number[7].

---

[3] Cf. [Ru53], p3, Definition 1.4.

[4] Cf. [Ru53], p39, Definition 3.10.

[5] See §4.1(*vii*) below.

[6] Cf. [Ru53], p39, Theorem 3.11.

Cantor's argument should, therefore, be interpreted only as establishing that we cannot define a number-theoretic function $f(n) = C_n$ in such a way that the range of $f(n)$ contains all, and only, the Cantorian real numbers.

However, there may yet be some definable number-theoretic function, $f(n)$, such that:

(*a*) its range contains all the Dedekind real numbers;

(*b*) it is undefined in the class R of Dedekind real numbers for some values in its domain;

(*c*) it is defined for these values in the extended class of Cantorian real numbers (which clearly contains all the Dedekind real numbers).

The conclusion: We can neither conclude that Cantor's diagonal argument determines an uncountable Dedekind real number, nor conclude from it that the cardinality of the Dedekind real numbers necessarily differs from that of the Dedekind (Peano) natural numbers[8].

## 3. Gödel's proof

In the second case, Gödel's use of the diagonal argument, in his seminal 1931 paper on formally undecidable propositions [Go31a], is purely illustrative. He uses it in his Introduction simply to sketch the main ideas of his Theorem VI heuristically, without claiming any rigour. The aim of his paper, as he then remarks, is to show that:

---

[7] Cf. [Ru53], p9, Theorem 1.32.

[8] We note that the latter result can be proved non-constructively in set theory (see [Ru53], p34, Theorem 2.40, Corollary). In ([An02], Corollary 1.1), however, we show that such proofs, which appeal to a Separation (Comprehension) Axiom (Schema) for introducing mathematical objects, may invite inconsistency.

"... The method of proof which has just been explained can obviously be applied to every formal system which, first, possesses sufficient means of expression when interpreted according to its meaning to define the concepts (especially the concept 'provable formula') occurring in the above argument; and, secondly, in which every provable formula is true. In the precise execution of the above proof, which now follows, we shall have the task (among others) of replacing the second of the assumptions just mentioned by a purely formal and much weaker assumption."

The question arises: Is it valid to treat Gödel's heuristic argument as being equivalent to his Theorem VI? Such equivalence is, of course, intuitionistically objectionable; even classically, set-theoretic arguments[9] corresponding to Cantor's diagonal argument are not considered constructive in any sense. That Gödel did not consider the two as equivalent is suggested by his remarks at the end of his Theorem VI:

"One can easily convince oneself that the proof we have just given is constructive (for all the existential assertions occurring in the proof rest upon Theorem V which, as it is easy to see, is intuitionistically unobjectionable), ...".

Thus, Gödel's conversion of the non-constructive, semantic, argument into a syntactic one that he considered constructive, and intuitionistically unobjectionable, is significant; although a strong Platonist, Gödel's remarks indicate that he intended to ensure that his argument, in Theorem VI, does not admit possible inconsistencies that could be argued as being inherent in the classically accepted non-constructivity of Cantor's diagonal argument.

---

[9] Cf. [Ru53], p34, Theorem 2.40, Corollary.

## 4.1 Turing real numbers

In the third case, the real numbers considered by Turing, in his original 1936 paper on computable numbers [Tu36], were explicitly defined as the (non-terminating) computations of "circle-free a-machines". Turing wrote:

(*i*) The "computable" numbers may be described briefly as the real numbers whose expressions as a decimal are calculable by finite means.

(*ii*) According to my definition, a number is computable if its decimal can be written down by a machine.

(*iii*) The computable numbers do not, however, include all definable numbers, and an example is given of a definable number which is not computable.

(*iv*) If at each stage the motion of a machine ... is *completely* determined by the configuration, we shall call the machine an "automatic machine" (or *a*-machine).

(*v*) If an *a*-machine prints two kinds of symbols, of which the first kind (called figures) consists entirely of 0 and 1 (the others being called symbols of the second kind), then the machine will be called a computing machine. If the machine is supplied with a blank tape and set in motion, starting from the correct initial *m*-configuration, the subsequence of the symbols printed by it which are of the first kind will be called the *sequence computed by the machine*. The real number whose expression as a binary decimal is obtained by prefacing this sequence by a decimal point is called the *number computed by the machine*.

(*vi*) If a computing machine never writes down more than a finite number of symbols of the first kind it will be called *circular*. Otherwise it is said to be *circle-free*.

(*vii*) A sequence is said to be computable if it can be computed by a circle-free machine. A number is computable if it differs by an integer from the number computed by a circle-free machine.

(*viii*) We shall avoid confusion by speaking more often of computable sequences than of computable numbers.

(*ix*) A computable sequence O is determined by a description of a machine which computes O. Thus the sequence 001011011101111... is determined by the table on p.234, and, in fact, any computable sequence is capable of being described in terms of such a table.

Again, as in the case of Cantorian real numbers, we cannot assume without a formal proof that a given Turing real number, defined as above, is also a Dedekind real number.

## 4.2 Turing's Halting argument

Although Turing appears to argue in his 1936 paper that Cantor's argument can be taken to establish the Platonic existence of an uncomputable Turing real number[10], he seems to have been ambivalent about using the argument unrestrictedly whilst introducing the Halting problem. Perhaps, echoing Gödel's reservations to some extent, Turing too felt the need to express - albeit obliquely - awareness of a non-constructive element in the use of Cantor's diagonal argument to define a (Dedekind) real number. Whatever the reason, he offered an alternative argument that was, essentially, based on defining an uncomputable number-theoretic function, rather than on non-constructively postulating that a period, followed by a non-terminating sequence of the digits 0 and 1, necessarily defines a Dedekind real number.

---

[10] Strictly speaking, Turing's Halting argument only implies that, whereas every Turing real number is, obviously, a Cantorian real number, some Cantorian real numbers are not Turing real numbers.

The ambivalence is reflected below, in Turing's 1936 description of the Halting problem:

(*x*) This new description of the machine may be called the *standard description* (S.D). It is made up entirely from the letters "*A*", "*C*", "*D*", "*L*", "*R*", "*N*", and from "**;**".

If finally we replace "*A*" by "1", "*C*" by "2", "*D*" by "3", "*L*" by "4", "*R*" by "5", "*N*" by "6", and "**;**" by "7" we shall have a description of the machine in the form of an arabic numeral.

(*xi*) The integer represented by this numeral may be called a *description number* (D.N) of the machine. The D.N determine the S.D and the structure of the machine uniquely. The machine whose D.N is *n* may be described as M(*n*).

(*xii*) To each computable sequence there corresponds at least one description number, while to no description number does there correspond more than one computable sequence. The computable sequences and numbers arc therefore enumerable.

(*xiii*) A number which is a description number of a circle-free machine will be called *satisfactory* number. In §8 it is shown that there can be no general process for determining whether a given number is satisfactory or not.

(*xiv*) ... we might apply the diagonal process. "If the computable sequences are enumerable, let In be the *n*-th computable sequence, and let Yn(*m*) be the *m*-th figure in In. Let J be the sequence with 1 − Yn(*n*) as its *n*-th figure. Since J is computable, there exists a number *K* such that 1 − Yn(*n*) = YK(*n*) all *n*. Putting *n* = *K*, we have 1 = 2YK(*K*), *i.e.* 1 is even. This is impossible. The computable sequences are therefore not enumerable".

The fallacy in this argument lies in the assumption that J is computable. It would be true if we could enumerate the computable sequences by finite means, but the problem of enumerating computable sequences is equivalent to the problem of finding out whether a given number is the D.N of a circle-free machine, and we have no general process for doing this in a finite number of steps. In fact, by applying the diagonal process argument correctly, we can show that there cannot be any such general process.

The simplest and most direct proof of this is by showing that, if this general process exists, then there is a machine which computes J. This proof, although perfectly sound, has the disadvantage that it may leave the reader with a feeling that "there must be something wrong". The proof which I shall give has not this disadvantage, and gives a certain insight into the significance of the idea "circle-free". It depends not on constructing J, but on constructing J ', whose $n$-th figure is Yn($n$).

($xv$) Let us suppose that there is such a process; that is to say, that we can invent a machine D which, when supplied with the S.D of any computing machine M will test this S.D and if M is circular will mark the S.D with the symbol "$u$" and if it is circle-free will mark it with "$s$. By combining the machines D and I we could construct a machine M to compute the sequence J'. The machine D may require a tape. We may suppose that it uses the $E$-squares beyond all symbols on $F$- squares, and that when it has reached its verdict all the rough work done by D is erased.

The machine H has its motion divided into sections. In the first $N-1$ sections, among other things, the integers 1, 2, …, $N-1$ have been written down and tested by the machine D. A certain number, say $R(N-1)$, of them have been found to be the D.N's of circle-free machines. In the $N$-th section the machine D tests the number $N$. If $N$ is satisfactory, $i.e.$, if it is the D.N of a circle-free machine, then $R(N) = 1 + R(N-1)$ and the first. $R(N)$ figures of the sequence of which a D.N is $N$ are calculated. The $R(N)$-th figure of this sequence is written down as one of the figures of the sequence  J' computed by H. If $N$ is not satisfactory, then $R(N) = R(N-1)$ and the machine goes on to the $(N+1)$-th section of its motion.

(*xvi*) From the construction of H we can see that H is circle-free. Each section of the motion of H comes to an end after a finite number of steps. For, by our assumption about D, the decision as to whether $N$ is satisfactory is reached in a finite number of steps. If $N$ is not satisfactory, then the $N$-th section is finished. If $N$ is satisfactory, this means that the machine M($N$) whose D.N is $N$ is circle-free, and therefore its $R(N)$-th figure can be calculated in a finite number of steps. When this figure has been calculated and written down as the $R(N)$-th figure of J', the $N$-th section is finished. Hence H is circle-free.

(*xvii*) Now let $K$ be the D.N of H. What does H do in the $K$-th section of its motion? It must test whether $K$ is satisfactory, giving a verdict "$s$" or "$u$". Since $K$ is the D.N of Hand since H is circle-free, the verdict cannot be "$u$". On the other hand the verdict cannot be "$s$". For if it were, then in the $K$-th section of its motion H would be bound to compute the first $R(K-1)+1 = R(K)$ figures of the sequence computed by the machine with $K$ as its D.N and to write down the $R(K)$-th as a figure of the sequence computed by H. The computation of the first $R(K) - 1$ figures would be carried out all right, but the instructions for calculating the $R(K)$-th would amount to "calculate the first $R(K)$ figures computed by $H$ and write down the $R(K)$-th". This $R(K)$-th figure wonld never be found. *I.e.*, H is circular, contrary both to what we have found in the last paragraph and to the verdict "$s$". Thus both verdicts are impossible and we conclude that there can be no machine D.

(*xviii*) The expression "there is a general process for determining …" has been used throughout this section as equivalent to "there is a machine which will determine …" This usage can be justified if and only if we can justify our definition of "computable". For each of these "general process" problems can be expressed as a problem concerning a general process for determining whether a given integer $n$ has a property $G(n)$ [*e.g.* $G(n)$ might mean "$n$ is satisfactory" or "$n$ is the Gödel 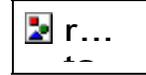 representation of a provable formula"], and this is equivalent to computing a number whose $n$-th figure is 1 if $G(n)$ is true and 0 if it is false.

(*xix*) The computable numbers do not include all (in the ordinary sense) definable numbers. Let P be a sequence whose $n$-th figure is 1 or 0 according as $n$ is or is not satisfactory. It is an immediate consequence of the theorem of §8 that P is not computable. It is (so far as we know at present) possible that any assigned number of figures of P can be calculated, but not by a uniform process. When sufficiently many figures of P have been calculated, an essentially new method is necessary in order to obtain more figures.

(*xx*) We cannot define general computable functions of a real variable, since there is no general method of describing a real number, but we can define a computable function of a computable variable.

## 4.3 The Halting problem and CT

Now, Turing's alternative argument is, essentially, the following form of the Halting problem[11]:

Suppose we have a Turing machine ***WillHalt*** which, given an input string ***M***+$w$, will halt and accept the string if Turing machine ***M*** halts on input $w$, and will halt and reject the string if Turing machine ***M*** does not halt on input $w$. Viewed as a Boolean function, **WillHalt**(**M**, *w*) (halts and) returns true in the first case, and (halts and) returns false in the second.

---

[11] Cf. David Matuszek's lecture notes at http://www.netaxs.com/people/nerp/automata/halting2.html.

Now, without appealing to Cantor's diagonal argument, we can prove by contradiction that Turing machine **WillHalt**(**M**, *w*) cannot exist. However, the proof does not address an implicit, and foundationally fundamental, non-constructive issue: what, precisely, do we mean by "**M** does not halt on input *w*"? If "**WillHalt**(**M**, *w*)" could effectively determine somehow that "**M** does not halt on input *w*", then this procedure could be built into **M**, and so **M** would halt on *w* (possibly with the perplexing annotation: "Sorry, this program is being halted as the program has determined that it does not halt!"). Thus, we cannot even define the machine **WillHalt**(**M**, *w*) as a Turing machine without inviting an immediate inconsistency!

Prima facie, the definition appears similar to that of a "***Liar***" proposition as "The '***Liar***' proposition is a lie", or that of a "***Russell***" set as "{*x* | '*x*' is a member of the '***Russell***' set if and only if '*x*' is not a member of '*x*'}"![12]

---

[12] We note that definitions such as those of the "***Liar***" proposition and of the "***Russell***" set are striking examples of what may be termed as pseudo-propositions of a language; they have no communicable content. They highlight the importance of viewing every formal system non-Platonically as a language of communication. Reasonably, any concept expressible within a communicable language must, ipso facto, be capable of effectively unambiguous interpretation.

Clearly, if the raison d'etre of a language is to be a medium for communicating content that lies outside the language, then the self-referential study of the internal structure of a language must also aim at effectively unambiguous communication - in this case, of the nature and structure of formal, finite, strings that are finitely constructible within the language.

Admitting non-finite, non-constructive elements into the language may compromise it to the point where we could permit non-communicable (or, more accurately, non-discernible) entities to be reflected in the language in such a way that they can be interpreted ambiguously by various interpretations - and yet give the illusion of unambiguous interpretation within each interpretation.

For instance, if we accept the argument that Cantorian real numbers are not necessarily Dedekind real numbers, then Cantor's diagonal construction may, indeed, be a notable example of such ambiguous interpretation.

So, if we accept that ***WillHalt***(***M***, *w*) is a well-defined Boolean function, then the significant conclusion from the above (and Turing's) argument is not that ***WillHalt***(***M***, *w*) is Turing-uncomputable, since we cannot define a Turing machine ***WillHalt***, but that, if the function is effectively computable, then the Church-Turing Thesis is false.

## 4.4 CT and Turing's "oracle"

We note that, in his 1936 paper, Turing did not conclude from his Halting argument that a function such as ***WillHalt***(***M***, *w*) is not computable; he only concluded that such a function cannot be computed by a Logical Computing (Turing) Machine as defined by him. However, since he was of the view that effective computability was, indeed, equivalent to computability by an LCM, Turing later proposed that a function such as ***WillHalt*** could, perhaps, be defined as an "oracle" that is not necessarily mechanical.

As noted by Hodges [Ho00]:

> Turing's 1938 Princeton Ph.D. thesis, work conducted in close cooperation with Church, was entitled *Systems of logic defined by ordinals,* and published as (Turing 1939). Predominantly the work consisted of highly technical developments within mathematical logic. However the driving force lay in the question: what is the consequence of supplementing a formal system with uncomputable deductive steps? In pursuit of this question, Turing introduced the definition of an "oracle" which can supply on demand the answer to the halting problem for every Turing machine. ... Turing defined the "oracle" purely mathematically as an uncomputable function, and said, "We shall not go any further into the nature of this oracle apart from saying that it cannot be a machine." The essential point of the oracle is that it performs *non-mechanical* steps.

Clearly, if the Church-Turing Thesis is false, then **WillHalt** can be a machine; further, if **WillHalt** is effectively computable, then CT is false. We thus consider whether there are mechanical computational processes that are not obviously duplicatable by the operations of a classical Turing machine.

Now, we can design a mechanical computer that will recognise a "looping" situation; it simply records every instantaneous tape description at the execution of each machine instruction, and compares the current instantaneous tape description with the record. We could then instruct the machine to assign arbitrary values to those undefined instances of **WillHalt**($M$, $w$) whose occurrences cause the machine to loop as in (*xvii*) above[13].

Is such a machine a classical Turing machine? Is the resultant predicate well-defined? Does the above construction transform a Turing-undecidable predicate into a decidable one? It is not obvious whether classical theory can provide unequivocal answers to these questions without recourse to CT, which would beg the question.

Clearly, if we can show, or accept, that the extended **WillHalt**($M$, $w$) is well-defined[14], then CT is false.

## 5. The influence of Gödel's reasoning

We can argue that the Platonist roots of Turing's "oracle" are dimly visible in the classical acceptance of the non-constructivity that is implicit in Gödel's reasoning. Prima facie, most classical interpretations uncritically draw sustenance, and strength, from the original reasoning and conclusions expressed by Gödel in his 1931 paper. Thus it would not be surprising if classically accepted heuristic, and formal, expositions are founded on implicit, non-constructive, Platonist assumptions that could be implicit in Gödel's formal, and informal, reasoning. The non-constructive logical consequences of such arguments would also, then, be classically accepted - perhaps reluctantly - simply because their non-constructive roots are obscurely hidden in such implicit premises.

---

[13] We note that, in the $K$'th section, machine H will input, and attempt to simulate, machine $K$ on input $K$ ad infinitum.

[14] That **WillHalt**($M$, $W$) is a well-defined function is classically well-accepted; the classical conclusion that the Halting argument constructively establishes the existence of a Turing-uncomputable number-theoretic function presumes such well-definedness.

For instance, we note that Gödel's and Rosser's reasoning is classically viewed as being definitively formal. The numerous attempts to challenge the reasoning and its conclusions over the years (starting with Wittgenstein most notably, and apparently led by Yessenin-Volpin later) have invariably, and quite effectively, attracted the fatal criticism of being heuristic and non-formal by Gödel's standard of constructiveness.

However, both Gödel's and Rosser's arguments appeal to semantic arguments at some point in their reasoning[15]. This could be obscured if, having accepted these results as definitively formal, classical expositions of these arguments appeal at some critical points to heuristic interpretations and reasoning for ease of exposition. Classical theory may, then, be focusing more on the implications of the accepted interpretations of Gödel's reasoning, and less on seeking any possible non-formal components that may be implicit either in the reasoning, or in its interpretation[16].

---

[15] See [An02a]. We argue there that, in Rosser's case, his application of the Deduction Theorem is invalid.

[16] For instance, see ([An02b], §4.4), and ([An02b], §4.5).